\documentclass[12pt,a4paper]{amsart}
\usepackage{amsfonts}
\usepackage{amsthm}
\usepackage{amsmath}
\usepackage{tikz}
\usepackage{amscd}
\usepackage[latin2]{inputenc}
\usepackage{t1enc}
\usepackage[mathscr]{eucal}
\usepackage{indentfirst}
\usepackage{graphicx}
\usepackage{graphics}
\usepackage{pict2e}
\usepackage[T1]{fontenc}
\usepackage{currvita}
\usepackage{url}

\usepackage{epic}
\numberwithin{equation}{section}
\usepackage[margin=2.9cm]{geometry}
\usepackage{epstopdf}

\usepackage{tikz,fullpage}
\usetikzlibrary{arrows,%
                petri,%
                topaths}%
\usepackage{tkz-berge}
\usepackage[position=top]{subfig}

\theoremstyle{plain}

 \theoremstyle{definition}

\newtheorem{?}[Th]{Problem}

\usetikzlibrary{graphs}
\usetikzlibrary{shadows}
\usetikzlibrary{graphs.standard}

\usetikzlibrary{graphs,graphs.standard}
\usetikzlibrary{positioning,chains,fit,shapes,calc}

\makeatletter
\usetikzlibrary{graphs,graphs.standard}
\newcount\nodecount
\tikzgraphsset{
  declare={subgraph N}%
  {
    [/utils/exec={\global\nodecount=0}]
    \foreach \nodetext in \tikzgraphV
    {  [/utils/exec={\global\advance\nodecount by1}, 
      parse/.expand once={\the\nodecount/\nodetext}] }
  },
  declare={subgraph C}%
  {
    [cycle, /utils/exec={\global\nodecount=0}]
    \foreach \nodetext in \tikzgraphV
    {  [/utils/exec={\global\advance\nodecount by1}, 
      parse/.expand once={\the\nodecount/\nodetext}] }
  }
}
\pgfmathdeclarefunction{alpha}{1}{%
  \pgfmathint@{#1}%
  \edef\pgfmathresult{\pgffor@alpha{\pgfmathresult}}%
}

\begin{document}
\author{Vojt\u{e}ch Dvo\u{r}\'ak}
\address{Trinity College, Cambridge CB21TQ, UK. Email address: vd273@cam.ac.uk}

\title[$P_{n}$-induced-saturated graphs exist for all $n \geq 6$]{$P_{n}$-induced-saturated graphs \\ exist for all $n \geq 6$
}


\begin{abstract} 
Let $P_{n}$ be a path graph on $n$ vertices. We say that a graph $G$ is $P_{n}$-induced-saturated if $G$ contains no induced copy of $P_{n}$, but deleting any edge of $G$ as well as adding to $G$ any edge of $G^{c}$ creates such a copy. Martin and Smith (2012) showed that there is no $P_{4}$-induced-saturated graph. On the other hand, there trivially exist $P_{n}$-induced-saturated graphs for $n=2,3$. Axenovich and Csik\'{o}s (2019) ask for which integers $n \geq 5$ do there exist $P_{n}$-induced-saturated graphs. R\"{a}ty (2019) constructed such a graph for $n=6$, and Cho, Choi and Park (2019) later constructed such graphs for all $n=3k$ for $k \geq 2$. We show by a different construction that $P_{n}$-induced-saturated graphs exist for all $n \geq 6$, leaving only the case $n=5$ open. 
\end{abstract}

\maketitle

\section{Introduction} 
Given graphs $G,H$, we say $G$ is \textit{$H$-saturated} if $G$ contains no subgraph isomorphic to $H$, but adding any edge from $G^{c}$ to $G$ creates a subgraph isomorphic to $H$. Related problems have been extensively studied (see for instance a survey of Faudree, Faudree and Schmitt \cite{faudree}).

In 2019, Axenovich and Csik\'{o}s \cite{axenovich} introduced the notion of induced-saturated graphs (before that in 2012, Martin and Smith \cite{martin} introduced a similar, more general notion). Given graphs $G,H$, we say $G$ is \textit{$H$-induced-saturated} if $G$ contains no induced subgraph isomorphic to $H$, but deleting any edge of $G$ creates an induced subgraph isomorphic to $H$, and adding any new edge to $G$ from $G^{c}$ also creates an induced subgraph isomorphic to $H$. Throughout the rest of the note, we will abbreviate a $H$-induced-saturated graph as a $H$-IS graph. 

While for any graph $H$, there exist $H$-saturated graphs, the same is not true for $H$-IS graphs. Indeed, for instance for a path on $4$ vertices $P_{4}$, Martin and Smith \cite{martin} showed that there exists no $P_{4}$-IS graph.

On the other hand, it is easy to see that there do exist $P_{2}$-IS and $P_{3}$-IS graphs. This leads to a question, asked by Axenovich and Csik\'{o}s \cite{axenovich}, for what integers $n \geq 5$ do there exist $P_{n}$-IS graphs. R\"{a}ty \cite{raty} was the first to make a progress on this question, showing by an algebraic construction that there exists a $P_{6}$-IS graph. Cho, Choi and Park \cite{cho} later showed that in fact for any $k \geq 2$, there exists a $P_{3k}$-IS graph. We use a different construction to settle the question completely, with the exception of the case $n=5$.

\vspace{5mm}

\textbf{Theorem 1.} For each $n \geq 6$, there is a $P_{n}$-induced-saturated graph.

\vspace{5mm}

In Section 2, we describe our construction of a $P_{n}$-IS graph $G_{n}$ for each $n \geq 6$. Then in Section 3, we check that the graph $G_{n}$ is actually $P_{n}$-IS.

\section{Construction}
We will construct, for each $n \geq 6$, a $P_{n}$-IS graph $G_{n}$. Our construction has been inspired by the observation of Cho, Choi and Park \cite{cho} that the Petersen graph is $P_{6}$-IS. We let

\begin{equation*}
V(G_{n})=\{ v_{1},...,v_{n-1},w_{1},...,w_{n-1} \}
\end{equation*}

Further, the edge set $E(G_{n})$ of $G_{n}$ is defined as follows. For $1 \leq i,j \leq n-1$, we have $v_{i}w_{j} \in E(G_{n})$ if and only if $i=j$. For $1 \leq i,j \leq n-1$, we have $v_{i}v_{j} \in E(G_{n})$ if and only if  $i-j \equiv \pm 1 \mod n-1$. And finally for $1 \leq i,j \leq n-1$, we have $w_{i}w_{j} \in E(G_{n})$ if and only if $i \neq j$ and $i-j \not\equiv \pm 1 \mod n-1$.

Note that the graph $G_{6}$ is isomorphic to the Petersen graph. Labelled graph $G_{7}$ is illustrated in the Figure 1 below, and (unlabelled) graphs $G_{6}$, $G_{7}$, $G_{8}$ are illustrated in the Figure 2 below.
\vspace{5mm}

\begin{center}
\begin{tikzpicture}
[
roundnode/.style={circle, draw=black, fill=white, minimum size=2mm},
squarednode/.style={rectangle, draw=red!60, fill=red!5, very thick, minimum size=5mm},
]
    \node[roundnode]  (v1) at ( 0, 0) {{\scriptsize $v_1$}}; 
    \node[roundnode] (v2) at ( 0,-1) {\scriptsize $v_2$};
    \node[roundnode] (v3) at ( 0,-2) {\scriptsize $v_3$};
    \node[roundnode] (v4) at ( 0,-3) {\scriptsize $v_4$};
    \node[roundnode] (v5) at (0,-4) {\scriptsize $v_5$};
    \node[roundnode] (v6) at ( 0, -5) {\scriptsize $v_6$}; 
    \node[roundnode] (w1) at ( 3, 0) {\scriptsize $w_1$};
    \node[roundnode] (w2) at ( 3,-1) {\scriptsize $w_2$};
    \node[roundnode] (w3) at ( 3,-2) {\scriptsize $w_3$};
    \node[roundnode] (w4) at ( 3,-3) {\scriptsize $w_4$};
    \node[roundnode] (w5) at ( 3, -4) {\scriptsize $w_5$}; 
    \node[roundnode] (w6) at ( 3, -5) {\scriptsize $w_6$};
    \begin{scope}[every path/.style={-}]
       \draw (v1) -- (w1);
       \draw (v2) -- (w2);
       \draw (v3) -- (w3);
       \draw (v4) -- (w4);
       \draw (v5) -- (w5);
       \draw (v6) -- (w6);
       \draw (v1) -- (v2);
       \draw (v2) -- (v3);
       \draw (v3) -- (v4);
       \draw (v4) -- (v5);
       \draw (v5) -- (v6);
    \end{scope}  
 \draw[-] (v1.west) .. controls  +(left:7mm)  .. (v6.west);
  \draw[-] (w1.east) .. controls +(right:30mm)  .. (w5.east);
  \draw[-] (w2.east) .. controls +(right:30mm)  .. (w6.east);
   \draw[-] (w1.east) .. controls  +(-30:18mm)  .. (w4.east);
   \draw[-] (w2.east) .. controls  +(-30:18mm)  .. (w5.east);
   \draw[-] (w3.east) .. controls  +(-30:18mm)  .. (w6.east);
   \draw[-] (w1.east) .. controls  +(-70:10mm)  .. (w3.east);
   \draw[-] (w2.east) .. controls  +(-70:10mm)  .. (w4.east);
   \draw[-] (w3.east) .. controls  +(-70:10mm)  .. (w5.east);
   \draw[-] (w4.east) .. controls  +(-70:10mm)  .. (w6.east);
\end{tikzpicture}

\vspace{5mm}

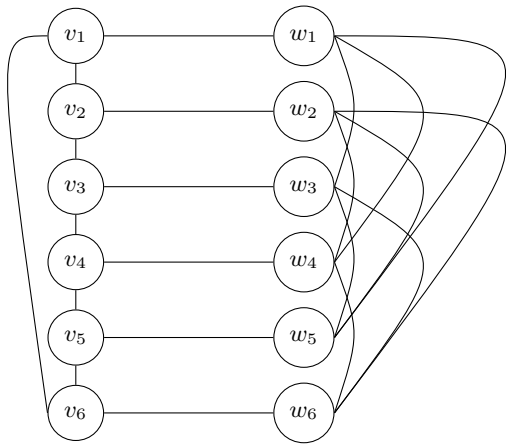
FIGURE 1. Labelled graph $G_{7}$
\vspace{5mm}

\tikzstyle{every node}=[circle, draw, fill=black!50,
                        inner sep=0pt, minimum width=4pt]
\begin{tikzpicture}[thick,scale=0.8,rotate=90]
    \draw \foreach \x in {0,72,...,288}
    {
        (\x:2) node {}  -- (\x+144:2)
        (\x:2) node {}  -- (\x+216:2)
        (\x:3) node {}  -- (\x+72:3)
        (\x:3) node {}  -- (\x:2)
    };
\end{tikzpicture}
\begin{tikzpicture}[thick,scale=0.8]
    \draw \foreach \x in {0,60,...,300}
    {
        (\x:2) node {}  -- (\x+120:2)
        (\x:2) node {}  -- (\x+180:2)
        (\x:3) node {}  -- (\x+60:3)
        (\x:3) node {}  -- (\x:2)
    };
\end{tikzpicture}
\begin{tikzpicture}[thick,scale=0.8,rotate=90]
    \draw \foreach \x in {0,360/7,720/7,1080/7,1440/7,1800/7,2160/7}
    {
        (\x:2) node {}  -- (\x+720/7:2)
        (\x:2) node {}  -- (\x+1080/7:2)
        (\x:3) node {}  -- (\x+360/7:3)
        (\x:3) node {}  -- (\x:2)
    };
\end{tikzpicture}
\vspace{5mm}

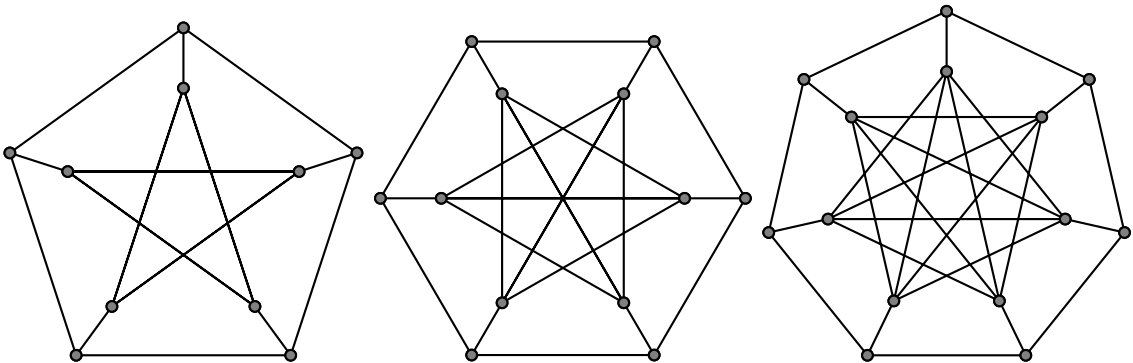
FIGURE 2. Graphs $G_{6}$, $G_{7}$, $G_{8}$
\vspace{5mm}
\end{center}

In the rest of the paper, we will prove that for each $n \geq 6$,  $G_{n}$ is $P_{n}$-IS, by checking the three properties that we need by the definition of an induced saturation.

\newpage

\section{Proof that the construction works}
\textbf{Claim 2.} For each $n \geq 6$, $G_{n}$ contains no induced copy of $P_{n}$.

\vspace{5mm}

\textit{Proof.} For $n=6$, the result is easy to check by hand. So throughout rest of the proof assume that $n \geq 7$. Also assume for contradiction that we have an induced copy of $P_{n}$ in $G_{n}$. 

First we claim that, since $n \geq 7$, among any five mutually disjoint vertices of the form $w_{i},w_{j},w_{k},w_{l},w_{m}$ for some $1 \leq i <j<k<l<m \leq n-1$, some three form a triangle $K_{3}$ in $G_{n}$. To see that, note that $G_{n}$ necessarily contains at least one of the edges $w_{i}w_{j},w_{j}w_{k},w_{k}w_{l},w_{l}w_{m},w_{m}w_{i}$ and due to the symmetry, we may without loss of generality assume that $G_{n}$ contains an edge $w_{i}w_{j}$. But then $w_{i}w_{j}w_{l}$ forms a triangle.

Write $W$ for $ \{ w_{1},...,w_{n-1} \} \subset V(G_{n})$. Since $P_{n}$ is acyclic, we must have at most four vertices from $W$ in our induced copy of $P_{n}$. We also must have at least one vertex from $W$ in our induced copy of $P_{n}$, since $|V(G_{n}) \setminus W|=n-1<n=|V(P_{n})|$. We will consider four cases depending on the number of vertices of $W$ in our induced copy of $P_{n}$.

If we have one vertex from $W$ in our induced copy of $P_{n}$, then we know our induced copy contains all of the vertices $v_{1},...,v_{n-1}$, but these form a cycle, which gives a contradiction.

If we have two vertices from $W$ in our induced copy of $P_{n}$, we may without loss of generality assume that our induced copy contains all of the vertices $v_{1},...,v_{n-2}$, but not the vertex $v_{n-1}$. Since $P_{n}$ contains no vertex of degree more than two, we know our copy of $P_{n}$ can not contain any of the vertices $w_{2},...,w_{n-3}$. But looking at all three two-element subsets of the set $\{ w_{n-2},w_{n-1},w_{1} \} $, we see that adding none of these subsets to the set $\{ v_{1},v_{2},...,v_{n-2} \} $ will create an induced copy of $P_{n}$.

Next assume we have three vertices $w_{i},w_{j},w_{k}$ from $W$ in our induced copy of $P_{n}$. Since $P_{n}$ is acyclic, we know we must have at least one of the relations $i-j \equiv \pm 1 \mod n-1, i-k \equiv \pm 1 \mod n-1, j-k \equiv \pm 1 \mod n-1  $ to hold, else $w_{i},w_{j},w_{k}$ would form a triangle. Consider two subcases depending on if one or two of the relations above hold (since $n>4$, we know all three can not hold simulateneously).

If two of the relations above hold, we may without loss of generality assume that we have precisely the vertices $w_{1},w_{2},w_{3}$ from $W$ in our induced copy of $P_{n}$. Then note that we must have $v_{2}$ in our copy of $P_{n}$ too, else the degree of $w_{2}$ in this copy would be zero. Also, $w_{2}$ has degree one in our copy of $P_{n}$, hence it forms one of the endpoints of $P_{n}$. Since $P_{n}$ is connected, our copy of it must also contain one of the vertices $v_{1}$ or $v_{3}$, due to the symmetry we may without loss of generality assume it contains $v_{1}$. But then it can not contain $v_{3}$, else it would contain a cycle $v_{1}v_{2}v_{3}w_{3}w_{1}$, hence $w_{3}$ also has degree one in our copy of $P_{n}$ and it forms another of the endpoints of $P_{n}$. But then we conclude $n \leq 5$, since distance of the endpoints of $P_{n}$ in our copy of it is at most four as we have a path $w_{3}w_{1}v_{1}v_{2}w_{2}$ connecting them, giving us a desired contradiction.

If just one of the relations above holds, we may without loss of generality assume that we have precisely the vertices $w_{1},w_{2},w_{j}$ for some $j$ such that $4 \leq j \leq n-2$ from $W$ in our induced copy of $P_{n}$. We can not have $v_{j}$ in our copy, else $w_{j}$ would have degree three in the copy. As $P_{n}$ is connected and $n>3$, we must have either $v_{1}$ or $v_{2}$ in our copy, and we can not have both, as then it would contain a cycle $v_{1}v_{2}w_{2}w_{j}w_{1}$. If we have $v_{1}$ but not $v_{2}$ in our copy of $P_{n}$, we can easily see that as $P_{n}$ is connected and $j \geq 4$, it can contain none of the vertices $v_{2},v_{3},...,v_{j}$, and hence contains at most $n-1$ vertices, giving us a contradiction. If we have $v_{2}$ but not $v_{1}$ in our copy, we conclude analogously by noting our copy of $P_{n}$ contains none of the vertices $v_{1},v_{n-1},...,v_{j}$ and $j \leq n-2$. 

Finally assume we have four vertices $w_{i},w_{j},w_{k},w_{l}$ from $W$ in our induced copy of $P_{n}$. If one of these four vertices is connected to all of the others, we must have a triangle in our copy of $P_{n}$ (since at least one of the three pairs of the other three vertices is connected too) and reach a contradiction. So due to this observation and the symmetry, it is enough to consider configurations $w_{1},w_{2},w_{l},w_{l+1}$ where $3 \leq l \leq n-2$. 

First consider the case $l=3$ (the case $l=n-2$ is analogous). In that case, we can not have $v_{1}$ or $v_{4}$ included in our copy of $P_{n}$, since that would mean degree of $w_{1}$ or $w_{4}$ respectively in the copy would be at least three. But then as $P_{n}$ is connected, none of the vertices $v_{4},v_{5},...,v_{n-2},v_{n-1},v_{1}$ can be in the copy, so our path $P_{n}$ has at most six vertices and hence $n \leq 6$, which is a contradiction.

Finally consider the case $3<l<n-2$. In this case $w_{1}w_{l}w_{2}w_{l+1}$ is a cycle, contradicting that $P_{n}$ is acyclic. $\square$

\vspace{5mm}

\textbf{Claim 3.} For each $n \geq 6$, deleting any edge of $G_{n}$ creates an induced copy of $P_{n}$.

\vspace{5mm}

\textit{Proof.} The edge we delete can be one of three types: $v_{i}v_{j}$, $v_{i}w_{j}$ or $w_{i}w_{j}$ for some $1 \leq i,j \leq n-1$; we consider these cases separately. 

First assume we delete an edge of the form $v_{i}v_{j}$. Then we must have $i-j \equiv \pm 1 \mod n-1$ and due to the symmetry, we may without loss of generality assume that the edge we deleted was $v_{1}v_{n-1}$. 

Then for $S_{1}= \{ w_{1} \} \cup \{ v_{i} : 1 \leq i \leq n-1 \}$, $G_{n}[S_{1}]$ is isomorphic to $P_{n}$.

Next assume we delete an edge of the form $v_{i}w_{j}$. Then we must have $i=j$, and due to the symmetry, we may without loss of generality assume that the edge we deleted was $v_{1}w_{1}$. 

Then for $S_{2}= \{ w_{1},w_{n-2} \} \cup \{ v_{i} : 1 \leq i \leq n-2 \}$, $G_{n}[S_{2}]$ is isomorphic to $P_{n}$.

Finally assume we delete an edge of the form $w_{i}w_{j}$ for some $i,j$ such that $i \neq j$, $i-j \not\equiv \pm 1 \mod n-1 $. Due to the symmetry, we may without loss of generality assume that the edge we deleted was $w_{1}w_{j}$ for some $j$ such that $3 \leq j \leq n-2$. 

Then if $3<j<n-2$, for $S_{3}= \{ w_{1},w_{j-1},w_{j},w_{n-1} \} \cup \{ v_{i} : 1 \leq i \leq j-2 \} \cup \{ v_{i} : j \leq i \leq n-3 \}$, $G_{n}[S_{3}]$ is isomorphic to $P_{n}$, if $j=3$, for $S_{3}'= \{ w_{1},w_{3} \} \cup \{ v_{1} \} \cup \{ v_{i} : 3 \leq i \leq n-1 \} $, $G_{n}[S_{3}']$ is isomorphic to $P_{n}$, and if $j=n-2$, for $S_{3}''= \{ w_{1},w_{n-2} \}\cup \{ v_{i} : 1 \leq i \leq n-2 \} $, $G_{n}[S_{3}'']$ is isomorphic to $P_{n}$.  $\square$

\vspace{5mm}

\textbf{Claim 4.} For each $n \geq 6$, adding any edge of $G_{n}^{c}$ to $G_{n}$ creates an induced copy of $P_{n}$.

\vspace{5mm}

\textit{Proof.} The edge we add can be one of three types: $v_{i}v_{j}$, $v_{i}w_{j}$ or $w_{i}w_{j}$ for some $1 \leq i,j \leq n-1$; we consider these cases separately. 

First assume we add an edge of the form $w_{i}w_{j}$. Then we must have $i-j \equiv \pm 1 \mod n-1$ and due to the symmetry, we may without loss of generality assume that the edge we added was $w_{1}w_{n-1}$. 

Then for $T_{1}= \{ w_{1},w_{n-1} \} \cup \{ v_{i} : 1 \leq i \leq n-2 \}$, $G_{n}[T_{1}]$ is isomorphic to $P_{n}$. 

Next assume we add an edge of the form $v_{i}v_{j}$ for some $i,j$ such that $i \neq j$, $i-j \not\equiv \pm 1 \mod n-1 $. Due to the symmetry, we may without loss of generality assume that the edge we added was $v_{1}v_{j}$ for some $j$ such that $3 \leq j \leq n-2$. 

Then if $3<j \leq n-2$, for $T_{2}= \{ w_{j-2},w_{j-1},w_{n-1} \} \cup \{ v_{i} : 1 \leq i \leq j-2 \} \cup \{ v_{i} : j \leq i \leq n-2 \}$, $G_{n}[T_{2}]$ is isomorphic to $P_{n}$, while if $j=3$, for $T_{2}'= \{ w_{2},w_{n-2},w_{n-1} \} \cup \{ v_{1}  \} \cup \{ v_{i} : 3 \leq i \leq n-2 \}$, $G_{n}[T_{2}']$ is isomorphic to $P_{n}$.  

Finally assume we add an edge of the form $v_{i}w_{j}$ for some $i \neq j$. Due to the symmetry, we may without loss of generality assume that the edge we added was $v_{1}w_{j}$ for some $j$ such that $2 \leq j \leq n-1$. 

Then if $2 \leq j \leq n-3$, for $T_{3}= \{ w_{j-1},w_{j},w_{j+1} \} \cup \{ v_{i} : 1 \leq i \leq j-1 \} \cup \{ v_{i} : j+1 \leq i \leq n-2 \}$, $G_{n}[T_{3}]$ is isomorphic to $P_{n}$, if $j=n-2$, for $T_{3}'= \{ w_{n-3},w_{n-2},w_{n-1} \} \cup \{ v_{i} : 1 \leq i \leq n-3 \} $, $G_{n}[T_{3}']$ is isomorphic to $P_{n}$ and if $j=n-1$, for $T_{3}''= \{ w_{n-2},w_{n-1} \} \cup \{ v_{i} : 1 \leq i \leq n-2 \} $, $G_{n}[T_{3}'']$ is isomorphic to $P_{n}$.  $\square$

\vspace{5mm}

This now completes the proof that the construction indeed works.

Finally, let us note that for $n=5$, our construction would give a graph that does contain an induced copy of $P_{5}$, leaving the question whether there exists a $P_{5}$-IS graph open.

\vspace{5mm}

\textbf{Note Added.} After the submission of this paper, it was pointed out to the author through personal correspondence that the case of $P_{5}$-IS graph was solved by Bonamy, Groenland, Johnston, Morrison and Scott in previously unpublished work. Their work can now be found in the following post \cite{bonamy}.

\section*{Acknowledgements}

The author would like to thank his PhD supervisor B\'ela Bollob\'as for his support and advice regarding the final version of this note.


\begin{thebibliography}{9}



\bibitem{axenovich}
Maria Axenovich, M\'{o}nika Csik\'{o}s.
\textit{Induced saturation of graphs}.
Discrete Mathematics \textbf{342}(4) (2018), 1195--1212.

\bibitem{bonamy}
Marthe Bonamy, Carla Groenland, Tom Johnston, Natasha Morrison, Alex Scott.
\textit{Induced Saturation for $P_{5}$}.
<https://tomjohnston.co.uk/blog/2020-05-22-induced-saturation-for-paths.html>.


\bibitem{cho}
Eun-Kyung Cho, Ilkyoo Choi, Boram Park.
\textit{On induced saturation for paths}.
arXiv:1907.05546 (2019).

\bibitem{faudree} 
Jill R. Faudree, Ralph J. Faudree, John R. Schmitt.
\textit{A Survey of Minimum Saturated Graphs}.
Electronic Journal of Combinatorics DS19 (2011) <http://eudml.org/doc/227276>.

\bibitem{martin}
Ryan R. Martin, Jason J. Smith.
\textit{Induced Saturation Number}.
Discrete Mathematics \textbf{312}(21) (2012), 3096--3106.

\bibitem{raty}
Eero R\"{a}ty.
\textit{Induced saturation of $P_{6}$}.
Discrete Mathematics \textbf{343}(1) (2020), 111641, 3 pp.





\end{thebibliography}
\end{document}